\documentclass[11pt,a4paper]{amsart}
\usepackage{amssymb,amsmath,amsthm}
\usepackage{graphicx}
\usepackage{enumerate}
\usepackage{tikz}

\newtheorem{theorem}{Theorem}[section]
\newtheorem{corollary}{Corollary}[section]
\newtheorem{lemma}{Lemma}[section]

\newtheorem{claim}{Claim}[section]

\newcommand{\C}{\mathcal{C}}
\newcommand{\D}{\mathcal{D}}
\newcommand{\R}{\mathbb{R}}
\newcommand{\h}{\mathcal{H}}

\newcommand{\al}{\alpha}

\newcommand{\lpos}{\mathrm{lpos} \,}
\newcommand{\pos}{\mathrm{pos} \,}
\newcommand{\pr}{\mathrm{pr} \,}
\newcommand{\lin}{\mathrm{lin} \,}

\begin{document}

\title{Positive bases, cones, Helly type theorems}

\author{Imre B\'ar\'any}

\begin{abstract} Assume that $k \le d$ is a positive integer and $\C$ is a finite collection of convex bodies in $\R^d$. We prove a Helly type theorem: If for every subfamily $\C^*\subset \C$ of size at most $\max \{d+1,2(d-k+1)\}$ the set $\bigcap \C^*$ contains a $k$-dimensional cone, then so does $\bigcap \C.$ One ingredient in the proof is another Helly type theorem about the dimension of lineality spaces of convex cones.
\end{abstract}


\maketitle

\section{Introduction and main result}\label{introd}

\bigskip
This paper is about Helly type properties of families of convex sets. Suppose for instance that $\C$ is a finite family of convex sets  in $\R^d$ and $\bigcap \C$ contains a halfline. Then of course $\bigcap \C^*$ also contains a halfline for every subfamily $\C^*$ of $\C.$ In the opposite direction assume that  $\bigcap \C^*$ contains a halfline for every subfamily $\C^*\subset \C$ of size at most $m.$ Does it follow that $\bigcap \C$ contains a halfline if $m=m(d)$ is chosen suitably? The answer is yes according to a theorem of Katchalski~\cite{Kat78}.

\begin{theorem}\label{th:ray} $m(d)=2d$, that is, if $\bigcap \C^*$ contains a halfline for every subfamily $\C^*\subset \C$ of size at most $2d$, then so does $\bigcap \C$.
\end{theorem}

A halfline is a one-dimensional cone. More generally a {\sl cone} $K$ with {\sl apex} $u \in \R^d$ and generator set $V\subset \R^d$ is the set of points of the form $u+\sum_{v\in W}^n\al(v)v$ where $W\subset V$ is finite and $\al(v)\ge 0$ for every $v \in W$. So $K$ consists of all finite and positive combinations of elements of $V$ translated by the vector $u$. The cone $K$ is {\sl polyhedral} if $V$ is finite. For properties of cones, their lineality spaces, their polar cones, etc. See for instance Gruber's book~\cite{gru} or Schneider's~\cite{schn}. 

\smallskip
Our first result extends Katchalski's theorem to $k$-dimensional cones where $k\in [d]=\{1,\ldots,d\}$. Set $m(k,d)=\max \{d+1,2(d-k+1)\}.$

\begin{theorem}\label{th:cone} Let $\C$ be a finite family of convex sets in $\R^d$ and $k\in [d]$. If $\bigcap \C^*$ contains a $k$-dimensional cone for every subfamily $\C^*\subset \C$ of size at most $m(k,d)$, then so does $\bigcap \C$.
\end{theorem}

Note that the case $k=0$ (which is not covered here) is Helly's theorem and we could have stated it by defining $m(0,d)=d+1.$

\medskip
The following examples show that the value of $m(k,d)$ is optimal. We write $ab$ for the scalar product of vectors $a,b\in \R^d.$

\medskip
{\bf Example 1.} Let $v_1,\ldots, v_{d+1}$ be the vertices of a regular simplex in $\R^d$ with $\sum_1^{d+1}v_i=0.$ Define $H_i$ as the halfplane $\{x \in \R^d: v_ix\le 0\},$ $i\in [d+1]$. Then $\bigcap_1^{d+1}H_i$ is a single point, namely the origin, so it contains no $k$-dimensional cone (for any $k>0$), but
for every $j\in [d+1]$ the set $\bigcap_{i\ne j}H_i$ contains a $k$-dimensional cone no matter what $k\in [d]$ is.

\medskip
{\bf Example 2.} Let $e_1,\ldots,e_d$ be the standard basis of $\R^d$ and define $H_i^+=\{x \in \R^d: e_ix\ge 0\}$ and $H_i^-=\{x \in \R^d: e_ix\le 0\}$ and set $\h_k=\{H_i^\pm: i\le d-k+1\}.$ Then $\bigcap \h_k$ is the subspace with $x_1=\ldots,x_{d-k+1}=0$, so it is a copy of $\R^{k-1}$ and can't contain a $k$-dimensional cone. Yet both $\bigcap (\h_k\setminus H_i^+)$ and $\bigcap (\h_k\setminus H_i^-)$ contain a $k$-dimensional cone (actually a $k$-dimensional halfspace) for every $i \in [d-k+1].$

\medskip
In these examples the convex sets in the family $\C$ are halfspaces of the form $\{x\in \R^d: ax\le 0\}.$ This is no coincidence: the proof of Theorem~\ref{th:cone} begins with a reduction from the family $\C$ to a family of halfspaces of this form. For the case of such halfspaces duality (or polarity) leads to a Helly type theorem on the dimension of the lineality space of cones (Theorem~\ref{th:pos}), which is the other main result in this paper.

\medskip
We remark that the long (and old and neat) survey paper by Danzer, Gr\"unbaum, and Klee~\cite{DGK} contains several similar Helly type results. The same applies to the more recent books by Gruber~\cite{gru}, Matou\v{s}ek~\cite{mat}, and Schneider~\cite{schn}. An up-to-date survey is by B\'ar\'any and Kalai~\cite{BarKal}.

\bigskip
\section{A Helly type theorem for the lineality space of cones}\label{sec:lineal}

Suppose $A \subset \R^d$ is a finite set and define $\pos A$ as the cone hull of $A=\{a_1,\ldots,a_n\}$, i.e., $\pos A=\{\sum_1^n\al_ia_i: \al_i\ge 0 \mbox{ for all }i\in [n]\}.$ The {\sl lineality space} of $\pos A$, to be denoted by $\lpos A$, is the set $\pos A \cap (-\pos A)$ which is a linear subspace of $\R^d$, actually the (unique) maximal dimension subspace that $\pos A$ contains, see for instance ~\cite{gru}. Set $h(k,d)=\max \{d+1,2(k+1)\}.$

\begin{theorem}\label{th:pos} Assume $A \subset \R^d$ is finite, $k \in [d]$, and $\dim \lpos B \le k$ for every $B\subset A$ with $|B|\le h(k,d).$ Then $\dim \lpos A \le k$ as well.
\end{theorem}

This is a Helly type result for the dimension of the lineality space. The value of $h(k,d)$ is optimal again as the following examples show. Define $A=\{v_1,\ldots,v_{d+1}\}$ with the $v_i$ coming from Example 1, then $\lpos A=\R^d$ while $\dim \lpos B=0$ for all $B\subset A$, $|B|\le d.$ This shows that $h(k,d)$ is optimal when $d+1\ge 2(k+1).$ For $d+1 \le 2(k+1)$ let $A_k$ be the set of vectors $\{\pm e_1,\ldots,\pm e_k\}.$ Now $\dim \lpos A_k=k$ but
$\dim \lpos A_k\setminus \{a\}<k$ for every $a \in A_k.$ These examples indicate a duality between Theorems~\ref{th:cone} and \ref{th:pos} that will become more transparent later.

\medskip
The proof method of Theorem~\ref{th:cone} works in other cases. For instance a theorem of Katchalski~\cite{Kat71} states the following.

\begin{theorem}\label{th:dim}  Let $\C$ be a finite family of convex sets in $\R^d$ and $k\in \{0,1,\ldots,d\}$. If $\dim \bigcap \C^* \ge k$ for every subfamily $\C^*\subset \C$ of size at most $m(k,d)$, then $\dim \bigcap \C \ge k$ as well.
\end{theorem}

The case $k=0$ is exactly Helly's theorem. We mention further that from this result Theorem~\ref{th:cone} (and in particular Theorem~\ref{th:ray}) can be deduced in a few lines. This is explained in Section~\ref{sec:app}.
Another example is the following result of de Santis \cite{deSan}.  

\begin{theorem}\label{th:deSan}  Let $\C$ be a finite family of convex sets in $\R^d$ and $k\in [d]$. If $\bigcap \C^*$ contains an affine flat of dimension $n-k$ for every subfamily $\C^*\subset \C$ of size at most $k+1$, then so does $\bigcap \C$.
\end{theorem}

\medskip
Examples 1 and 2 show that the bounds $m(k,d)$ (in Theorem~\ref{th:dim}) and $k+1$ (in Theorem~\ref{th:deSan}) are optimal. Section~\ref{sec:app} explains how the proof method of Theorem~\ref{th:cone} applies to the last two results.

\section{Proof of Theorem~\ref{th:pos}}

Define $L=\lpos A$ and choose a positive basis $X\subset A$ of $L$. This is just a subset of $A$ with $\pos X=L$ but $\pos (X\setminus \{x\})\ne L$ for any $x \in X.$ We are going to use a theorem of Reay~\cite{Reay} stating that a positive basis $X$ has a partition $X_1\cup \ldots \cup X_r$ such that $|X_1|\ge |X_2|\ge \ldots \ge |X_r|\ge 2$ and $X_1\cup \ldots \cup X_j$ is a positive basis for $\lin(X_1\cup \ldots \cup X_j)$ whose dimension is $|X_1\cup \ldots \cup X_j|-j$ for every $j \in [r].$ Note that this and $|X_i|\ge 2$ imply that $\dim \lin(X_1\cup \ldots \cup X_j)\ge j.$

\medskip
We show first that, in our case, $|X_1|\le k+1.$ As $X_1$ is a positive basis of $\lin X_1$ whose dimension is $|X_1|-1\le d$ we have $|X_1|\le d+1\le h(k,d).$ Then, according to the condition of the theorem, $|X_1|-1=\dim \pos X_1 \le k$ and $|X_1|\le k+1$ indeed.

\medskip
Set $B_j=X_1\cup \ldots \cup X_j$. We show, by induction on $j$, that $|B_j|\le k+j$ or that $\dim \lin B_j \le k.$ (The two are equivalent since $\dim \lin B_j=|B_j|-j$.)   The starting step $j=1$ was just fixed so we move to step $j-1 \to j$ assuming that $j\ge 2.$ Since $X_j$ has the smallest size among the sets $X_1, \ldots, X_j$ we have $|X_j|\le \frac 1{j-1}|B_{j-1}|$. By the induction hypothesis $|B_{j-1}|\le k+(j-1).$ Thus
\begin{eqnarray*}
j &\le& \dim \lin B_j=|B_j|-j=|B_{j-1}|+|X_j|-j\\
    &\le& k+(j-1)+\frac {k+(j-1)}{j-1}-j=\frac {jk}{j-1},
\end{eqnarray*}
and $j\le \frac {jk}{j-1}$  implying that $j-1\le k.$

\begin{claim} $\frac {jk}{j-1}+j \le h(k,d).$
\end{claim}

This claim implies that $|B_j|\le h(d,k)$ and then the condition of Theorem~\ref{th:pos} gives $\dim \lin B_j \le k$. Consequently $|B_j|\le k+j$, completing the induction.

\medskip
{\bf Proof.} Assume first that $h(k,d)=2(k+1)$ so $d+1\le 2(k+1).$ Direct computation shows that the inequality $\frac {jk}{j-1}+j \le h(k,d)$ is equivalent to $(j-1)(j-2)\le k(j-2)$ which is true for $j\ge 2$ because $j-1 \le k.$ Suppose next that $h(k,d)=d+1$, so $k\le \frac {d-1}2$. Then $\frac {jk}{j-1}+j \le d+1$ holds iff  $\frac {j(d-1)}{2(j-1)}+j \le d+1$. The last inequality is the same as $2j^2-5j+2 \le d(j-2)$, and here $2j^2-5j+2=(j-2)(2j-1)$ and $2j-1\le 2(k+1)-1\le d$ indeed.\qed

\medskip

For $j=r$ the Claim gives $\dim L= \dim \lin (B_r)=\dim \lpos A \le k$.\qed

\medskip
{\bf Remark.} A strange step in the proof is that the inequalities $|B_j|\le \frac {jk}{j-1}+j \le h(k,d)$ imply, via $ \dim \lin B_j=|B_j|-j$, that $|B_j|\le k+j$.

\medskip
\section{Proof of Theorem~\ref{th:cone}}

As $m(k,d)\le d+1$ Helly's theorem shows that $\bigcap \C$ is non-empty. We assume, after a translation if necessary, that $0 \in\bigcap \C$. We assume further that every $C \in \C$ is closed. This assumption is justified since a convex set contains a cone if and only if its closure contains a translate of this cone.

\begin{lemma}\label{l:apex} Let $C \subset \R^d$ be a closed convex set containing a cone $K$ with apex at $u.$ Then $C$ contains the cone $v-u+K$ for every $v \in C$; the apex of this cone is at $v$.
\end{lemma}

The {\bf proof} is simple and is omitted. We remark though that it is enough to check the case when $K$ is a halfline.\qed

\medskip
We agree that from now on the word ``cone'' means a cone with apex at the origin. In view of Lemma~\ref{l:apex}, $\C$ satisfies the condition

\medskip
\hskip 1cm (*)\hskip0.7cm  $\bigcap \C^*$ contains a $k$-dimensional cone for every $\C^* \subset \C$\\
\hspace*{3cm}    whose size is at most $m(k,d)$.

\medskip
After these preparations the proof starts by reducing or changing the family $\C$ in two steps. For the first step we choose a $k$-dimensional cone $K(\C^*)\subset \bigcap \C^*$ for every $\C^* \subset \C$ satisfying $|\C^*| \le m(k,d),$ and we choose $K(C^*)$ so that it has exactly $k$ (of course linearly independent) generators. Replace each $C \in \C$ by the cone hull (or what is the same, convex hull), $D(C)$, of the union of all cones $K(\C^*)$ with $C \in \C^*$ and set $\D=\{D(C): C\in \C\}.$ The new system $\D$ consists of polyhedral cones and satisfies condition (*). Moreover $\bigcap \D \subset \bigcap \C$ because $D(C) \subset C$. So it suffices to show that $\bigcap \D$ contains a $k$-dimensional cone.

\medskip
For the second reduction we observe that each $D(C)$ is the intersection of finitely many closed halfspaces, $H$, of the form $\{x\in \R^d: ax\le 0\}$ for some $a \in \R^d$ ($a \ne 0)$, the outer normal of $H.$ Let $\h$ be the collection of these (finitely many) closed halfspaces, and let $A \subset \R^d$ be the set of the corresponding outer normals. Evidently $\h$ satisfies condition (*) and $\bigcap \C$ contains a $k$-dimensional cone if so does $\bigcap \h=\bigcap \D$.

\medskip
The solution set of the system of linear inequalities
\[
ax\le 0, a\in A
\]
coincides with $\bigcap \h$ and then $\bigcap \h$ is the polar of the cone $K=\pos A.$ Let $L$ be the lineality space of $K$, then $K$ is the sum of $L$ and the cone $L^{\perp}\cap \pos A$, see~\cite{gru}. The latter cone is a pointed cone in $\R^d$ and in the subspace $L^{\perp}$ as well. It coincides with the cone hull of the orthogonal projection, to be denoted by $\pr A$, of $A$ onto $L^{\perp}.$ Using standard properties of polarity (c.f. \cite{gru}) we see that
\begin{eqnarray*}
\bigcap \h&=& K^{\circ}=L^{\perp}\cap (L^{\perp}\cap \pos A)^{\circ}\\
  &=&(L^{\perp}\cap (L+(\pos \pr A)^{\circ})=L^{\perp}\cap (\pos \pr A)^{\circ},
\end{eqnarray*}
where the polars are taken in $\R^d$. Thus $\bigcap \h$ is a cone in $L^{\perp}$ and is full dimensional there because $L^{\perp}\cap \pos A=\pos \pr A$ is a pointed cone. Then $\bigcap \h$ contains a $k$-dimensional cone if and only if $\dim L^{\perp} =\dim (\lpos A)^{\perp}$ is at least $k$.

\medskip
Condition (*) for $\h$ says that $\bigcap \h^*$ contains a $k$-dimensional cone for every $\h^* \subset \h$ with $|\h^*|\le m(k,d)$. Writing $B\subset A$ for the outer normals of the halfspaces in $\h^*$ we have $\bigcap \h^*=(\pos B)^{\circ}$. The previous argument applies to $B$ (instead of $A$) and gives that $\bigcap \h^*$ contains a $k$-dimensional cone if and only if $\dim (\lpos B)^{\perp} \ge k$.

\medskip
Then for every $B\subset A$ with $|B|\le m(k,d)$, we have  $\dim (\lpos B)^{\perp} \ge k$, or, what is the same  $\dim (\lpos B) \le d-k$. Theorem 3 applies now and shows that $\dim (\lpos A) \le d-k$, and equivalently $\dim L^{\perp}\ge k$. Thus $\bigcap \h$ contains a $k$-dimensional cone.\qed

\bigskip
A byproduct of this proof is an interesting corollary for systems of homogeneous inequalities.

\begin{corollary}  Assume $A$ is a finite set of nonzero vectors in $\R^d.$ The system $ax\le 0,\; a\in A$ has at least $k$ linearly independent solutions if and only if for every $B\subset A$ whose size is at most $m(k,d)$ the system $ax\le 0,\; a\in B$ has at least $k$ linearly independent solutions.
\end{corollary}

\section{Proofs of Theorems~\ref{th:dim} and \ref{th:deSan}}\label{sec:app}

{\bf Proof} of Theorem~\ref{th:dim}. We only give a sketch. Again each $C \in \C$ is closed and contains the origin. For every subfamily $\C^*$ of size at most $m(k,d)$ we choose a set of $k$ linearly independent vectors from $\bigcap \C^*$. They together with the origin show that $\bigcap \C^*$ is indeed at least $k$-dimensional.

\medskip
Next each $C \in \C$ is replaced by the cone hull, $D(C)$, of those $k$ tuples of vectors that were chosen for some $\C^*$ with $C \in \C^*.$  Let $\D$ be the family of all $D(C)$. This time the condition on $\C$ implies that for every $\D^* \subset \D$ of size at most $m(k,d)$, the set $\bigcap \D^*$ contains a $k$-dimensional cone. Then, according to Theorem~\ref{th:cone}, there is  a $k$-dimensional cone in $\bigcap \D$ implying that there are linearly independent vectors $v_1,\ldots,v_k \in \bigcap \D$. It is not hard to check (we omit the details) that for a small enough $r>0$ the vectors $rv_i \in \bigcap \C$ for all $i\in [k]$. So $\bigcap \C$ is also at least $k$-dimensional.\qed

\bigskip
We explain next how Katchalski's theorem (Theorem~\ref{th:dim}) implies Theorem~\ref{th:pos}. We assume again that every $C \in \C$ is closed and contains the origin. The same reduction as in the proof of Theorem~\ref{th:pos} works again: for each $\C^*\subset \C$ whose size is at most $m(k,d)$ we choose the same cone $K(C^*)\subset \bigcap \C^*$ whose generators are $k$ linearly independent vectors from $\bigcap \C^*.$ Replace each $C \in \C$ by the cone hull, $D(C)$, of the union of all cones $K(\C^*)$ with $C \in \C^*$ and set $\D=\{D(C): C\in \C\}.$ As we have seen $D(C)\subset C$ for every $C\in \C.$ Then $\dim \bigcap \D^*\ge k$ for every $\D^* \subset \D$ whose size is at most $m(k,d).$ Thus Theorem~\ref{th:dim} implies that $\dim \bigcap \D\ge k$ and then there are linearly independent vectors $v_1,\ldots,v_k \in \bigcap \D \subset \bigcap \C.$ The cone hull of these vectors is a cone of dimension at least $k$ in $\bigcap \C$ finishing the proof.

\bigskip
There is a related result by Ko{\l}odziejczyk~\cite{Kol} stating a condition, similar to (*), guaranteeing that $\bigcap \C$ contains an $k$-dimensional affine halfspace. The condition is that $\bigcap \C^*$ contains an $k$-dimensional affine halfspace for every $\C^* \subset \C$ whose size is at most $m(k,d)$. This theorem can also be proved by the same method.

\bigskip
{\bf Proof} of Theorem~\ref{th:deSan}. We assume again that each $C\in \C$ is closed and that the origin lies in every $C \in \C.$ Check that if $C$ contains an affine flat $F$, then it also contains the subspace $F-f$ where $f\in F$ is arbitrary. The condition of the theorem is then modified to the following: For every subfamily $\C^* \subset \C$ whose size is at most $k+1$, $\bigcap \C^*$ contains an $(n-k)$-dimensional subspace. We have to show that $\bigcap \C$ contains an $(n-k)$-dimensional subspace.

\medskip
Now comes the reduction in two steps. For each such $\C^*$ choose such a subspace $F(\C^*)$ and replace every $C \in \C$ by the convex hull, $D(C)$, of all subspaces $F(\C^*)$ with $C \in C^*.$ It is easy to check that every $D(C)$ is a subspace in $\R^d$. The new system $\D=\{D(C): C \in \C\}$ satisfies the previous condition, namely, that for every subfamily $\D^* \subset \D$ whose size is at most $k+1$, $\bigcap \D^*$ contains an $(n-k)$-dimensional subspace. It suffices to show that $\bigcap \D$ contains a $(n-k)$-dimensional subspace.

\medskip
Every $D \in \D$ is a subspace and one can choose $(1+\dim D)$ closed halfspaces whose intersection is $D$. Fix these halfspaces for every $D$ and let $\h$ be the collection of these halfspaces. So $H\in \h$ is of the form $\{x \in \R^d: ax\le 0\}$ where $a \in \R^d$ is the outer normal of the halfspace $H$. Write $A$ for the set of all outer normals in $\h.$

\medskip
As $\bigcap \h =\bigcap \D$, our target is to show that $\bigcap \h$ contains an $(n-k)$-dimensional subspace. The condition is that the intersection of any $k+1$ of these halfspaces contains an $(n-k)$-dimensional subspace. This happens if and only if the outer normals of these $k+1$ halfspaces are linearly dependent. So the condition says that there is no linearly independent $(k+1)$-element subset in $A$, in other words, $\dim \lin A \le k$ which implies in turn that $\bigcap \h$ contains an $(n-k)$-dimensional subspace.\qed

\bigskip
{\bf Acknowledgements.} This piece of work was partially supported by Hungarian National Research Grants No 131529, 131696, and 133819.

\bigskip
\vspace{.5cm} {\sc Imre B\'ar\'any}\\[1mm]
  {\footnotesize Alfr\'ed R\'enyi Institute of Mathematics}\\[-1mm]
  {\footnotesize 13-15 Re\'altanoda Street, Budapest, 1053 Hungary,}\\[-1mm]
 {\footnotesize e-mail: {\tt imbarany@gmail.com}, and }\\
  {\footnotesize Department of Mathematics, University College London}\\[-1mm]
  {\footnotesize Gower Street, London WC1E 6BT, UK}\\

\end{document}